

\baselineskip=14pt
\parskip=10pt

\font\eighttt=cmtt8
\magnification=\magstephalf

\def\1{{\overline{1}}}
\def\2{{\overline{2}}}
\parindent=0pt
\overfullrule=0in

\def\frac#1#2{{#1 \over #2}}
\bf
\centerline
{ Another Hanukkah Miracle:}
\centerline
{The Gaps Between Consecutive  Christmas-in-Hanukkah Years is ALWAYS a Fibonacci Number!
}
\rm
\bigskip
\centerline{ \it By Lisa BUDD and Shalosh B. EKHAD}

\qquad \qquad \qquad \qquad \qquad
{\it In fond memory of Marvin Isadore Knopp (1933-2011) z''l, who left us on Christmas Eve, the fifth candle of
Hanukkah 5772, who did {\bf not} like computers (even for Email), but was nevertheless one of the kindest and
warmest {\bf neshamot} (souls) on this planet.}

{\bf Preface}: The first version of this paper, authored only by the second-named author, left the question of the title
open. The first-named author found a very convincing explanation. Chapter I is the original version, and Chapter II is
the new material due to the first-named author.

{\bf Chapter I}

Jane Legrange, the beloved wife of my beloved master, commented this last Christmas Eve that it is so nice
to be able to light a {\it menorah} (or {\it chanukiah}) at the same time as the (Christian) {\it goyim} 
sitting next to their Christmas trees, and she asked her husband, who asked me, how special is that coincidence?

Using the Maple package 

{\tt http://www.math.rutgers.edu/\~{}zeilberg/tokhniot/LUACH }

written by Zeilberger a few years ago, and updated and upgraded recently in order to answer this question,
I figured out that in this, third, millenium, it happens $270$ times, so \%27 of the time. In the next two
millenia after that it would happen $266$ times in each. In the sixth millenium (i.e. between 5001 and 6000),
it would only happen $263$ times, while in the seventh millenium (i.e. between 6001 and 7000), it would only
happen $258$ times. If we are still alive in the eighth millenium (i.e. between 7001 and 8000), it would only
happen $134$ times, and in the ninth millenium, only 15 times! After that it will never happen for many many
years, certainly not until 20000 AD. The last year is 8473, where we would light the first candle on
Christmas Eve. After that, Hanukkah will always be way {\bf after Christmas}.

But cheer up! Eventually we can, every few years, sit in a Sukkah during Christmas Eve!
The first time this would happen is 16103, when Christmas would be on the 20th day of Tishrei,
and we can put the Christmas tree inside the Sukkah. During the 17th millenium (between 16001 and 17000)
we would only have $68$ Sukkah Christmases, but during the 18th (between 17001 and 18000) we would have
$207$ of them, during the 19th we would have $239$ but only $235$ during the 20th millenium.

While doing these calculations, I noticed something really amazing! Between the Gregorian years of 1801
and 7390 (when Christmas-in-Hanukkha is starting to get rarer and rarer until it disappears completely)
the gap between consecutive such lucky years is {\bf always}  a Fibonacci number, in fact, it is always
a member of the set $\{ 2,3,5,8\}$. It would be interesting to have a {\it conceptual} proof, in addition
to my computational proof.

And guess what! Once Christmas-in-Sukkot becomes not too rare, and that would happen in 17064 A.D., we have the
same phenomenon, the gaps between consecutive Christmas-in-Sukkot years is also {\bf always} a
member of the set $\{ 2,3,5,8\}$, at any rate, until 20000 A.D. Of course, eventually Christmas will
say good-bye to Sukkot also, and start visiting Rosh Hashana (but not that often since
the latter is only two days long), then Shavuot, and then Pesach.  Readers are welcome to 
experiment with the Maple package {\tt LUACH} to continue these preliminary investigations.

For a list of all the Christmas-in-Hanukkah years (between 1801 and 20000 (of course it ends in 8478)) see,

{\tt http://www.math.rutgers.edu/\~{}zeilberg/tokhniot/oLUACHs1 }

and for a list of all the Christmas-in-Sukkot years (also between 1801 (of course it only starts in 16103) and 20000), see

{\tt http://www.math.rutgers.edu/\~{}zeilberg/tokhniot/oLUACHs2 } .

{\bf Chapter II}

The answer to the question in the title is not to do with Fibonacci numbers, as shown by the fact that 1 is also a Fibonacci number.  The phenomenon
is due to the possible lengths in days of periods of 1, 2, 3, etc. years.  A Jewish year may have 353, 354, 355, 383,
384 or 385 days.  By contrast, a Gregorian year may have 365 or 366 days.  Years divisible by four have 366 days, except
for century years not divisible by 400 such as 1900 or 2100, and all other years have 365 days.

As the Jewish and Gregorian years are of different lengths, the date in the Gregorian calendar of a given Jewish date will vary from year to year.  Typically, after a 354-day ordinary Jewish year it will be 365-354 = 11 days earlier in the Gregorian calendar and after a 384-day Jewish leap year it will be 384-365 = 19 days later in the Gregorian calendar.  As a result, in the short term, a given Jewish date may during several decades fall on any of about 30 different Gregorian dates, and vice versa.  Further, in the long term Jewish dates are on average getting later compared with Gregorian dates.  This is because the average length of a Jewish year is about 365.246822 days, but the average for a Gregorian year is 365.2425 days.  This difference amounts to over six minutes a year or 4.3 days per thousand years.  After about 84,500 years, the difference amounts to a whole year, so the starts of periods during which there can be the concidence of Christmas and 
Hanukkah are separated by about 84,500 years.

At present, the Jewish festival of Hanukkah or Chanukah always falls close to Christmas.  It lasts for eight days, and between 1900 and 2099 the Gregorian date of the first day ranges from 28th November to 27th December.  Clearly, as Hanukkah lasts for eight days, one of these days coincides with Christmas if and only if the first day is 18th to 25th December inclusive.  In Chapter I, Ekhad calls years when this happens Christmas-in-Hanukkah years. It notes that during a period of several millennia, from 1801 to 7390, gaps between such years are always 2, 3, 5 or 8 years.  Actually, this is the case since the introduction of the Gregorian calendar in 1582 and, extrapolating that calendar backwards (and if the current rules for the Jewish calendar applied then), it would have been true since 876.  Of course, 2, 3, 5 and 8 are all in the Fibonacci sequence.  However, the analysis below shows that this is just a coincidence.  Before or after 876 to 7390, the calendar drift means that Christmas-in-Hanukkah years were or will be less frequent, so the gaps between them tend to be bigger.  The longer gaps in the periods shortly before and after 876 to 7390 are 11 and 19 years, not the next Fibonacci numbers, 13 and 21.

As noted above, Jewish years can have any of six lengths.  Further, the lengths of periods of Jewish years measured in days fall into two sets differing by about 30 days, depending on how many leap years the periods contain (unless the period is a multiple of 19 years).  These may be called the shorter Jewish periods and the longer Jewish periods.  The possible lengths in days for various periods of Jewish years are tabulated in 
{\tt http://hebrewcalendar.tripod.com/lengths2.html }.  We thus have the following list of lengths of periods in days (Shorter Jewish period, then Gregorian period, then Longer Jewish period):

{\bf 1 year:} 353-355; 365-366; 383-385
{\bf 2 years:} 707-710; 730-731; 737-740
{\bf 3 years:} 1092-1094; 1095-1096; 1121-1124
{\bf 4 years:} 1445-1449; 1460-1461; 1475-1477
{\bf 5 years:} 1799-1803; 1825-1827; 1829-1832
{\bf 6 years:} 2184-2187; 2190-2192; 2214-2217
{\bf 7 years:} 2538-2541; 2555-2557; 2569-2571
{\bf 8 years:} 2892-2895; 2921-2922; 2922-2925
{\bf 9 years:} 3276-3279; 3286-3288; 3306-3309
{\bf 10 years:} 3630-3633; 3651-3653; 3661-3664
{\bf 11 years:} 4015-4018; 4016-4018; 4044-4048
{\bf 12 years:} 4368-4372; 4382-4383; 4399-4402

{\bf Possible and impossible intervals}

Why cannot two consecutive years both be Christmas-in-Hanukkah years?  Consider periods of one year.  The lengths in the Jewish calendar must be at least 365-355 = 10 days shorter or 383-366 = 17 days longer than in the Gregorian one.  This means that if Christmas falls during Hanukkah in the first year, it must fall well before or well after it in the following year, so the years cannot both be Christmas-in-Hanukkah.

Similarly, gaps of four, seven, 10 and 12 years are impossible.  Gaps of two years are just about possible, as a longer Jewish period may be only six or seven days longer than a Gregorian one.

Gaps of six years seem possible, as a Gregorian period may be as little as three days longer than a shorter Jewish one.  However, if there is a Christmas-in-Hanukkah after six years, there must also be one half-way through the period, so this becomes two gaps of three years.  The same is true of gaps of nine years, as a Gregorian period may be as little as seven days longer than a shorter Jewish one, but they become three gaps of three years.

Finally, gaps of more than ten years are impossible during the peak period, as is shown by the fact that gaps of 11 years are possible going by lengths of periods, but in fact none occurs.

It follows that the only possible gaps are 2, 3, 5 and 8 years, as observed.

{\bf Festivals other than Hanukkah}

With the continual drift of the Jewish calendar relative to the Gregorian, Christmas will eventually fall close to each Jewish festival in turn.  The same analysis as above applies to Pesach (Passover), which is also eight days.  Succot (Tabernacles) is nine days, but the calculations by ``Shalosh B. Ekhad", confirmed by an analysis similar to the above, show that this makes little difference.  For a two-day festival, such as Rosh Hashanah (New Year) or Shavuot (Pentecost), gaps between consecutive Christmas-in-festival years must necessarily be greater.  They cannot be less than eight years; they may be 11, 19 or more years even at peak times.

\bigskip
\hrule

\bigskip

Lisa Budd, London, England,

Email: {\eighttt lisabudd at gmail dot com} ,

\bigskip
Shalosh B. Ekhad, c/o D. Zeilberger, Department of Mathematics, Rutgers University (New Brunswick),
Hill Center-Busch Campus, 110 Frelinghuysen Rd., Piscataway,
NJ 08854-8019, USA.

Email:  {\eighttt ShaloshBEkhad at gmail dot com} ,

Written: Jan. 2, 2012.  Revised: Nov. 2, 2019.

\end